\documentclass[twocolumn]{cinc}

\usepackage{graphicx}
\usepackage{amsmath}
\usepackage{amssymb}
\usepackage{url}
\usepackage{xcolor}
\usepackage{wrapfig}

\def\v{{\mathbf{v}}}
\def\w{{\mathbf{w}}}
\def\g{{\mathbf{g}}}
\def\u{{\mathbf{u}}}
\def\ds{\mathrm{d}{\it s}}
\def\dx{\mathrm{d}{\it x}}
\def\I{\textbf{I}}
\def\F{\textbf{F}}
\def\Diff{\rm D}
\begin{document}
\bibliographystyle{cinc}

% Keep the title short enough to fit on a single line if possible.
% Don't end it with a full stop (period).  Don't use ALL CAPS.
\title{Hyperelastic Modeling of Cardiac Tissue Deformation for Data Registration}

% Both authors and affiliations go in the \author{ ... } block.
% List initials and surnames of authors, no full stops (periods),
%  titles, or degrees.
% Don't use ALL CAPS, and don't use ``and'' before the name of the
%  last author.
% Leave an empty line between authors and affiliations.
% List affiliations, city, [state or province,] country only
%  (no street addresses or postcodes).
% If there are multiple affiliations, use superscript numerals to associate
%  each author with his or her affiliations, as in the example below.

\author {Joyce Ghantous$^{1}$, Yves Coudière$^{1}$, Michael Leguèbe$^{1}$ \\
\ \\ % leave an empty line between authors and affiliation
 $^1$ Univ. Bordeaux, CNRS, Inria, Bordeaux INP, IMB, UMR 5251, IHU Liryc, F-33400 Talence, France }

\maketitle

% LaTeX inserts the ``Abstract'' heading in the proper style and
% sets the text of the abstract in italics as required.

\begin{abstract}
This work proposes a biomechanical registration framework designed to combine functional optical-mapping da\-ta with structural high-resolution MRI of cardiac tissue. The two acquisitions are performed in different mechanical configurations: optical mapping is obtained after stretching the sample, whereas MRI is acquired in a fixed resting state. To estimate the deformation between these configurations, the tissue is modeled as a nearly incompressible transversely isotropic hyperelastic material, and the traction amplitude is identified from optical-map marker positions through an inverse problem. Two- and three-dimensional simulations show physically admissible deformations and provide a first proof of concept for the proposed strategy.
\end{abstract}

% LaTeX inserts the extra space here automatically.

\section{Introduction}

Optical mapping and high-resolution MRI provide complementary information on cardiac tissue. Optical mapping records surface electrical activity often in a mechanically stretched configuration. High-resolution MRI acquired at 9.4T is performed with the tissue in a resting configuration. It provides detailed three-dimensional microstructural information. This prevents a direct connection between functional and anatomical data. This problem differs from standard multimodal registration problems, such as MRI--CT registration, where both images generally describe the same three-dimensional anatomy, despite differences in image contrast. Classical image-registration methods, including diffeomorphic or statistical-shape approaches such as~\cite{Bone2018Deformetrica}, and biomechanical approaches such as equilibrated warping~\cite{Genet2018EquilibratedWarping} relying on image-based information, are not suited for the present setting.

\begin{figure}[htbp!]
    \centering
     \includegraphics[height=0.48\linewidth]{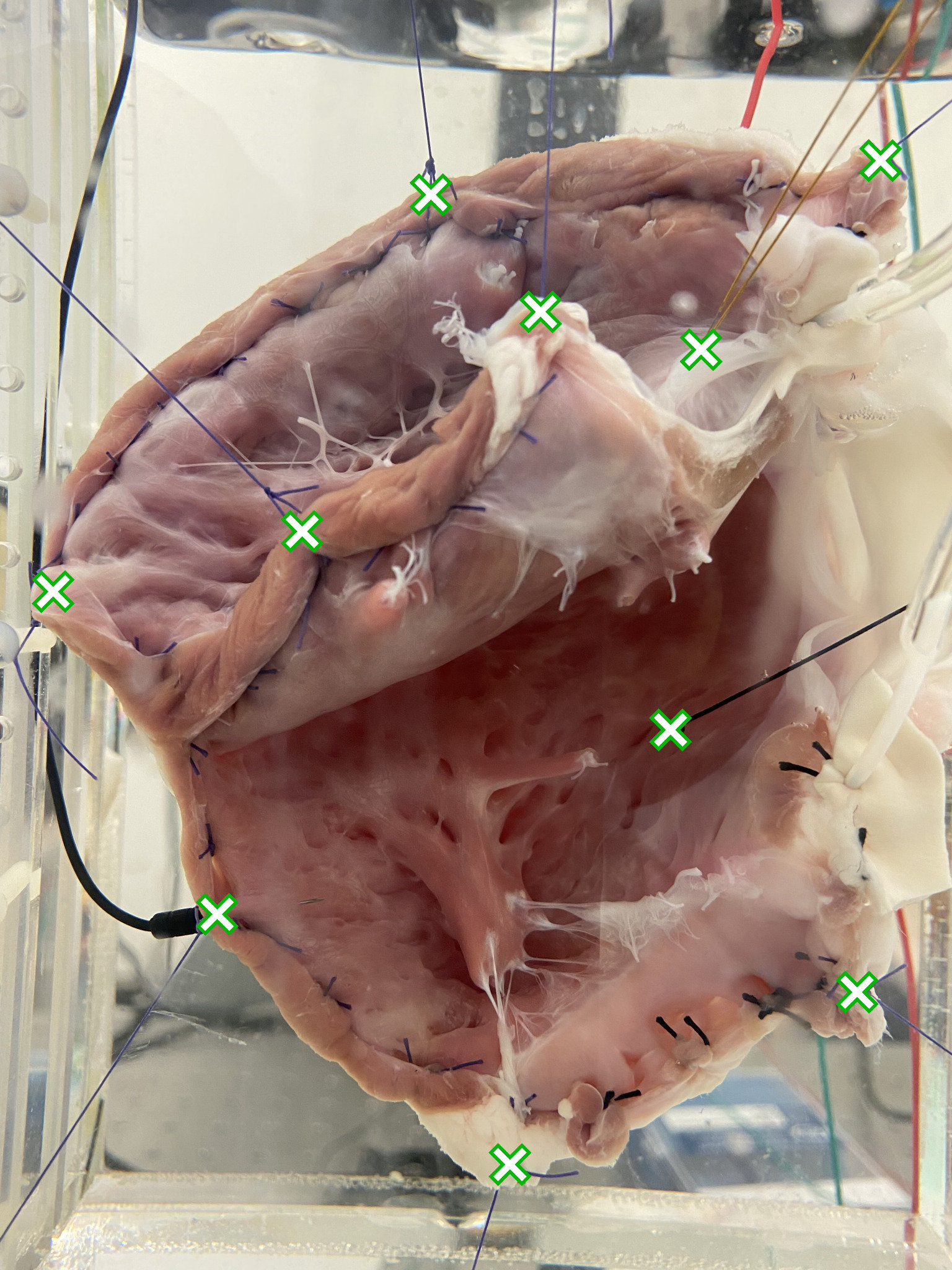} \quad
    \includegraphics[height=0.48\linewidth]{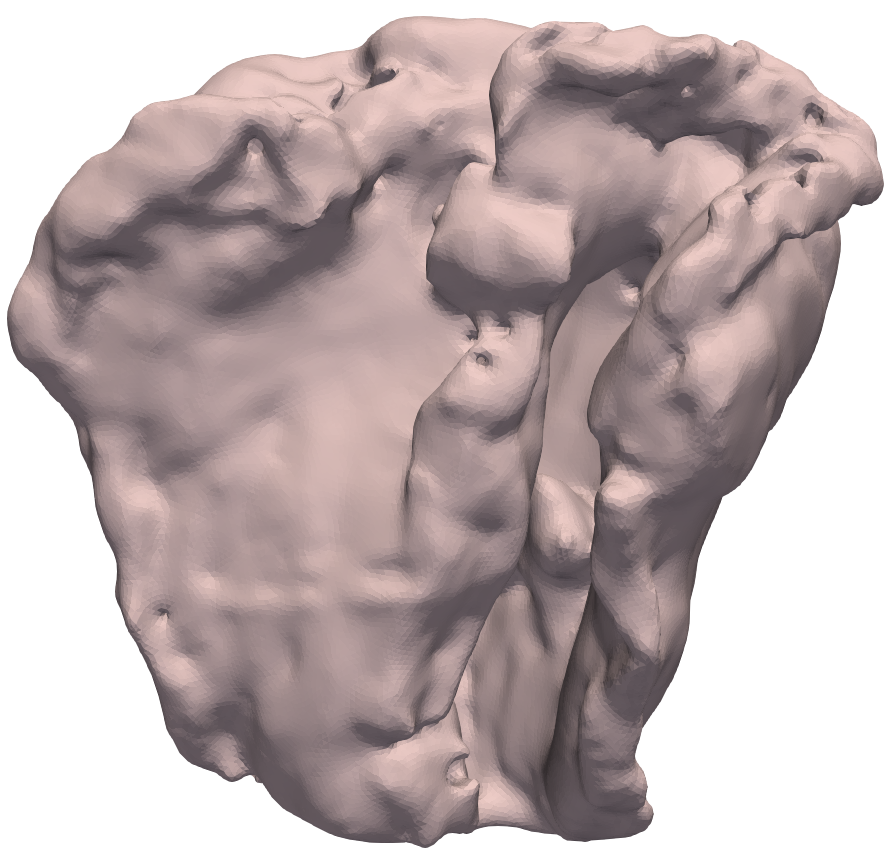}
    \caption{\textcolor{black}{Image of the stretched wedge preparation with attachment points indicated by crosses (left), and a 3D display of the same sample in its resting state (right).}}
    \label{fig:3d-rv-mesh}
\end{figure}

In this work, we propose a biomechanical registration strategy to map a resting configuration obtained with an MRI setup onto a stretched optical-mapping configuration. A key application of this strategy is to combine data acquired from sheep cardiac tissue wedge preparations composed of the right ventricle and part of the septum, as illustrated in Figure \ref{fig:3d-rv-mesh} and detailed in  \cite{walton2018compartmentalized,pannetier2023modeling}. The tissue is modeled as a nearly incompressible transversely isotropic hyperelastic material, following cardiac tissue models \cite{HolzapfelOgden2009,gasser2006hyperelastic}. Boundary tractions are applied to mimic the wires used during the optical-mapping experiment to stretch the sample. The resulting large deformation is computed with a finite element code, and the traction amplitude is estimated through an inverse problem based on selected marker positions. Preliminary 2D and 3D simulations provide a first proof of concept before applying the method to the full experimental geometries (Figure \ref{fig:3d-rv-mesh}). 

The purpose of this model is to generate a physically admissible deformation consistent with the experimental stretching and suitable for data registration. Ultimately, this framework would support the validation of the pacemaker model in \cite{pannetier2023modeling}, while providing a useful tool for experimentalists to investigate mechanisms related to arrhythmias.

% \textbf{Mechanical model.} Cardiac tissue has a complex microstructure, where cardiomyocytes are embedded in a collagen- and elastin-rich extracellular matrix. This organization leads to a nonlinear, anisotropic, and nearly incompressible mechanical behavior, usually described within finite-strain continuum mechanics. Passive myocardial mechanics is often modeled using structure-based hyperelastic laws, such as the Holzapfel--Ogden model \cite{HolzapfelOgden2009}, which takes into account the fiber direction. This formulation is used in ventricular simulations, residual-stress modeling, and inverse parameter identification \cite{wang2013structure,wang2014modified,Gao2015}. More recent works have also investigated myocardial compressibility under multiaxial loading, motivating nonlinear volumetric extensions of anisotropic hyperelastic models \cite{Avazmohammadi2020compressibility}.

\section{Materials and method}
\label{sec:materials}
\textbf{Kinematics and constitutive model.} Let $\Omega_0 \subset {\mathbb{R}}^d$, $d=2,3$, be the reference domain representing an MRI geometry. Its boundary is decomposed into three disjoint parts $\Gamma=\Gamma_D\cup\Gamma_N\cup\Gamma_0$, with no displacement on $\Gamma_D$, prescribed traction on $\Gamma_N$, and traction-free conditions on~$\Gamma_0$. We denote the deformation map $\varphi:\Omega_0 \to {\mathbb{R}}^d$:
\[
\varphi(x) := x + \u(x), \qquad x \in \Omega_0,
\]
with $\u:\Omega_0 \to {\mathbb{R}}^d$ denoting the displacement field. 
The deformation gradient is $\F := \nabla \varphi = \I + \nabla \u$, and its Jacobian $J := \det \F$ measures the local volume change.

The cardiac tissue is modeled as nearly incompressible and transversely isotropic, more rigid along the muscle fibers. We use the hyperelastic energy density (see~\cite{HolzapfelOgden2009,gasser2006hyperelastic}):
\begin{multline*}
    % W (\tilde I_1, \tilde I_2, \tilde I_4, J)
    \textcolor{black}{W(\tilde I_1,\tilde I_{d-1},\tilde I_4,J)}: = \sum_{p=1}^{d-1} \mu_p(\tilde I_p-d)
+ C_0\,e^{C_1(\tilde I_1-d)^2} \\
+ C_2\,e^{C_3(\tilde I_4-1)^2}
+ \kappa(J-1)-\kappa\ln(J),
\end{multline*}
where $\mu_p, C_0, C_1, C_2, C_3,\kappa>0$ are material parameters reported in Table \ref{tab:param}. The reduced invariants $\tilde{I}_p$, $p=1,\dots,d-1$, describe the isochoric isotropic deformation of the tissue matrix, whereas $\tilde{I}_4$ captures the stretch along the fiber direction. 
Their full definitions are given in \cite[Eq. 3.5--3.6]{HolzapfelOgden2009} and \cite[\S 3.1]{gasser2006hyperelastic}. 

\textbf{The forward problem.} Let $\g:\Gamma_N\to{\mathbb{R}}^d$ be a traction on a non-empty boundary part $\Gamma_N\subset\Gamma$. In the hyperelastic setting, we seek a minimizer of the total energy as follows
\begin{equation*}
    \u \in \arg\min_{\w\in V^+} \mathrm{E}(\w),
\end{equation*}
where the admissible set is defined by %$V^{+}:=\left\{\u\in [{\rm H}^1(\Omega_0)]^d:\ J>0 \,;\ \u= 0\ \text{on }\Gamma_D\right\},$ 
$V^{+}:=\{\w\in [{\rm H}^1(\Omega_0)]^d:\ J(\w)>0,\ \w=0 \text{ on }\Gamma_D\},$
and the total energy $\mathrm{E}:V^{+}\to{\mathbb{R}}$ is given by
${\mathrm{E}}(\u) := \int_{\Omega_0} W(\u) \,\dx - \int_{\Gamma_N}\g\cdot\u\,\ds.$

The nonlinearity of the hyperelastic constitutive law leads to a nonlinear forward problem. It is formulated by seeking a stationary point $\u\in V^+$ of the total energy  such that
\begin{equation}
    \label{eq:diff_problem}
\Diff\mathrm{E}(\u;\v)=0
\qquad \forall\,  \v \in [H^1_{0,\Gamma_D}(\Omega_0)]^d,
\end{equation}
where $\Diff\mathrm{E}(\u;\v)$ denotes the differential of $\mathrm{E}$ at $\u$ in the direction $\v$ and is given by $\Diff \mathrm{E}(\u;\v)
=
\int_{\Omega_0} \mathbf P(\u):\nabla\v\,\dx
-
\int_{\Gamma_N}\g\cdot\v\,\ds,$ where $\mathbf P(\u)$ is the first Piola--Kirchhoff stress tensor.

\textbf{The inverse problem.} For a prescribed traction vector $\g$, the
hyperelastic model \eqref{eq:diff_problem} generates a deformed domain. In this work, the
traction is assumed to act along the $x$-direction, namely
$\g = g\,\mathbf{e}_x$, where $\mathbf{e}_x$ denotes the unit vector in the~$x$-direction and the scalar amplitude $g$ lies in the admissible interval
$[g_{\min}, g_{\max}] = [0,4]$~N/cm$^2$.
The choice of bounds for $g$ is motivated by the experimental setup: strong enough to stretch the tissue, but not too large to avoid damage.
%the lower bound corresponds to the minimum traction needed to stretch the tissue sample, whereas the upper bound is chosen small enough to avoid damaging the tissue. % or causing it to detach from the optical-mapping frame. 
% The values are not directly measured, but inferred from the optical mapping markers.

In order to track deformation, markers are placed in the sample and visible in both modalities (e.g., toothpicks; see Figure \ref{fig:optical-map}). \textcolor{black}{We denote the coordinate of the $i$-th marker by $x_i$ in the deformed configuration (more precisely their projection of the pixel grid). We search for the traction amplitude $g$ which minimizes the discrepancy between the simulated position $x_i(g)$ and the measured ones ${x_i}^*$.}

In the 2D case (Figure~\ref{fig:2d-inverse-sketch}), we generate a set of synthetic target coordinates % $\{{x_i}^*\}$ by computing a deformation with a prescribed force $g_{\mathrm{ref}}$: 
%Only the coordinates of some markers via toothpicks in the optical map are available. In the 2D case in Figure~\ref{fig:2d-inverse-sketch}, we use their $x$-coordinates and we denote $x_i^*$ the target $x$-coordinate of the $i$-th marker. In the synthetic case, these target coordinates are generated by perturbing the resulting forward-model positions obtained for a reference amplitude $g_{\mathrm{ref}}$, namely
$x_i^* = x_i(g_{{\mathrm{ref}}}) + \varepsilon_i$, by computing a deformation with a prescribed force $g_{\mathrm{ref}}$ and a uniform measurement uncertainty $
\varepsilon_i \sim {U} \left(
-\frac{L}{2N_{\mathrm{pix}}},
\frac{L}{2N_{\mathrm{pix}}}
\right),
$, where $L$ is the length of the optical map domain along the $x$-axis and $N_{{\mathrm{pix}}}=256$ is the number of pixels in the interval.

We look for the force $g^*$ which minimizes the following weighted least-squares functional:
\begin{equation}
    \label{eq:costFct}
     F(g) = \sum_{i=1}^{N_x} w_i \bigl(x_i^* - x_i(g)\bigr)^2,
\quad w_i = \frac{|x_i^*|^q}{\sum_{j=1}^{N_x} |x_j^*|^q}.
\end{equation}
%where $x_i(g)$ is a simulated $x$-coordinate of the $i$th marker for a given amplitude $g$.
The exponent $q\geq 0$ controls the relative weight $w_i$ assigned to each marker; in this work, we set $q=1$. This choice assigns larger weights to markers with larger target coordinates. 

\textbf{Numerical Framework.}
The reference domain $\Omega_0 \subset {\mathbb{R}}^d$ is discretized by a mesh composed of triangles in 2D and tetrahedra in 3D. We solve Problem~\eqref{eq:diff_problem} with ${\mathbb{P}}^1$-Lagrange finite elements. 

The material parameters (Table \ref{tab:param}) are chosen from values reported for passive myocardial tissue, and are adjusted to ensure plausible deformations (see \cite{Guccione1991,HolzapfelOgden2009,gasser2006hyperelastic,ogden1978nearly}). The parameters $\mu_p$, $C_0$, and $C_1$ control the isotropic contribution of the tissue matrix, whereas $C_2$ and $C_3$ weight the anisotropic response associated with the local fiber direction. The parameter $\kappa$ penalizes the volume variation imposing the quasi-incompressibility constraint. 
\begin{table}[htbp!]
\centering
\renewcommand{\arraystretch}{0.9}
\begin{tabular}{lccccccc}
\hline
 & $\mu_1$ & $\mu_2$ & $C_0$ & $C_1$ & $C_2$ & $C_3$ & $\kappa$ \\
\hline
2D case & 1   & --  & 0.1809 & 9   & 0.058 & 15  & 80 \\
3D case & 1 & 1 & 0.1809 & 9   & 0.058 & 15  & 80 \\
% 3D case & 0.1 & 0.1 & 10   & 0.5 & 10    & 0.9 & 40  \\
\hline
\end{tabular}
\caption{Material parameters used in the numerical simulations, in $\mathrm{N/cm^2}$ for $\mu_1$, $\mu_2$, $C_0$, $C_2$, and $\kappa$; $C_1$ and $C_3$ are dimensionless.}
\label{tab:param}
\end{table}

The nonlinear forward problem is solved using a Newton--Raphson method combined with a Goldstein line-search strategy. At each Newton iteration, the tangent linear system is assembled and solved using the parallel direct solver MUMPS. The Newton iterations are stopped when the energy gradient norm is smaller than the chosen tolerance ($10^{-8}$). The biomechanical code is implemented in FreeFem++\footnote{\url{https://freefem.org/}} using MPI together with the ffddm library. % The mesh is partitioned with Scotch, and the forward 2D simulations are performed on a standard desktop computer using 6 MPI processes in approximately 2 $\mathrm{min}$, whereas the 3D forward simulation requires about $1\mathrm{h}\,  57 \, \mathrm{min}$ on $32$ MPI processes on a cluster.

The inverse problem is solved in Python using the L-BFGS-b algorithm suitable for bound-constrained optimization problems. At each iteration, the current trial value $g$ is passed to the forward solver through a Python wrapper. The forward hyperelasticity problem is then solved, and the resulting marker positions are compared with the target coordinates $\{x^*_i\}$. The corresponding value of~$F(g)$ is returned to the optimizer until convergence. %The 3D and inverse problem computations are performed using the MCIA (Mésocentre de Calcul Intensif Aquitaine) computing facilities.

\section{Results}

%This section presents the numerical results obtained for the forward and inverse problems in 2D with an extension to the 3D forward problem. 

\begin{figure}[htbp!]
    \centering
    \includegraphics[width=0.73\linewidth]{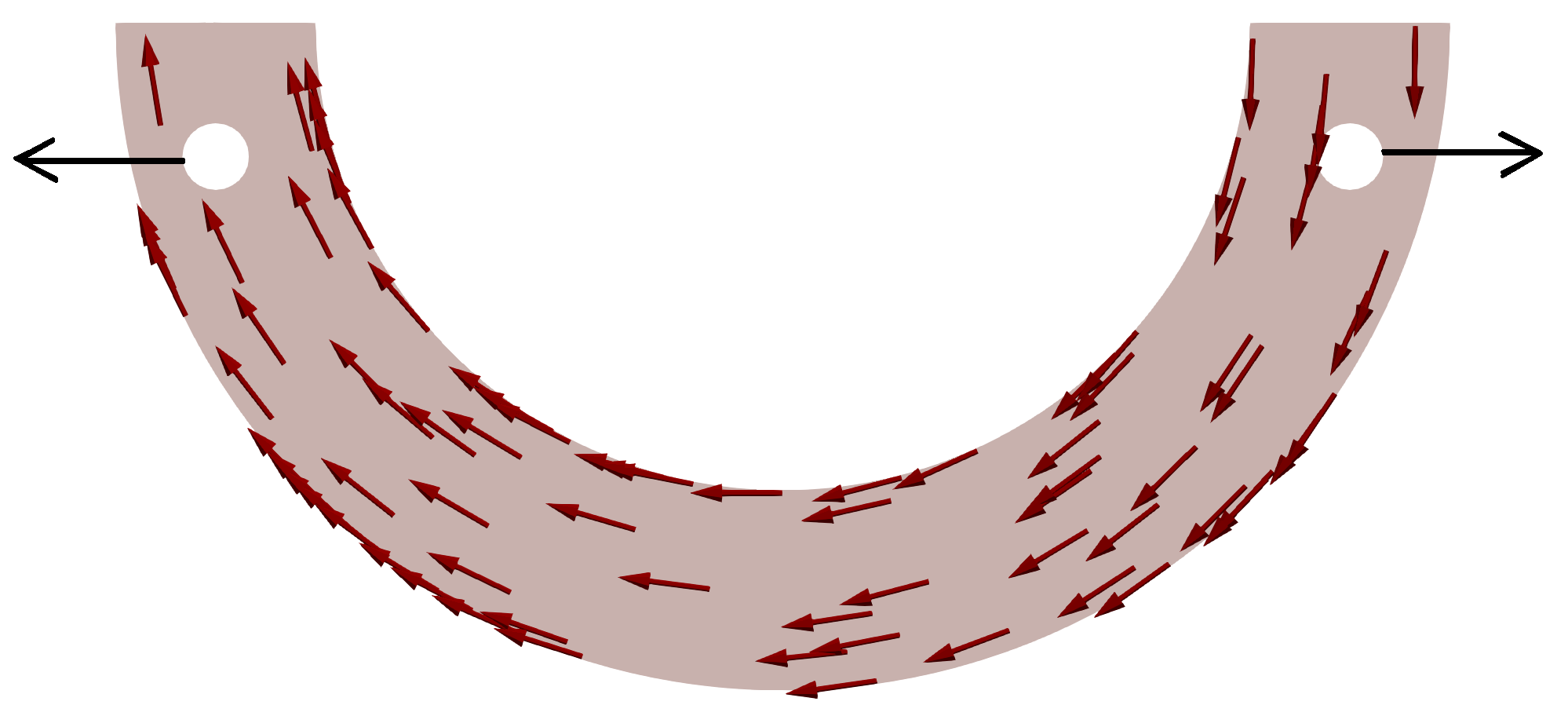}

    {\large $\downarrow$}
    \vspace{-0.05cm}

    \includegraphics[width=0.87\linewidth]{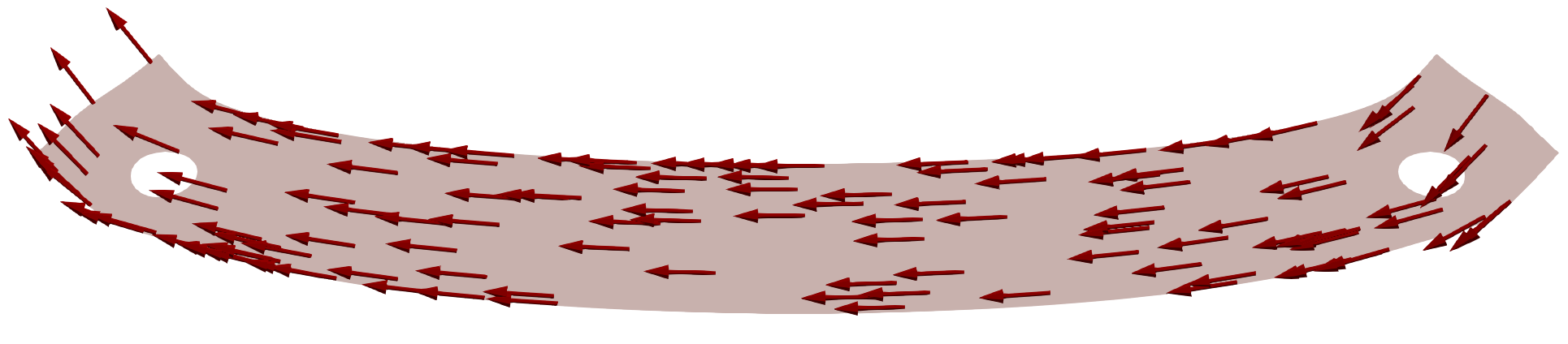}

   \caption{2D forward problem: reference domain with vectors representing the horizontal traction at the top, and deformed domain. \textcolor{black}{A sample of the reference and deformed local fiber directions is shown in red for both cases}.}
    \label{fig:2d-forward-deformation}
\end{figure}

\textbf{Forward problem on a 2D domain.} The considered domain consists of a
lower semi-annular domain subjected to horizontal traction, as illustrated in Figure~\ref{fig:2d-forward-deformation}. These forces are applied on the inner disks representing the wires passed through the tissue sample during the optical mapping. Applying the displacement field obtained by solving the forward problem \eqref{eq:diff_problem} leads to a flattened configuration in Figure~\ref{fig:2d-forward-deformation}. The relative volume variation remains small ($0.97\%$), which is consistent with the quasi-incompressible behavior commonly assumed for soft biological tissues~\cite{Fung1993Biomechanics}.

\begin{figure}[htbp!]
    \centering
    \includegraphics[height=0.45\linewidth]{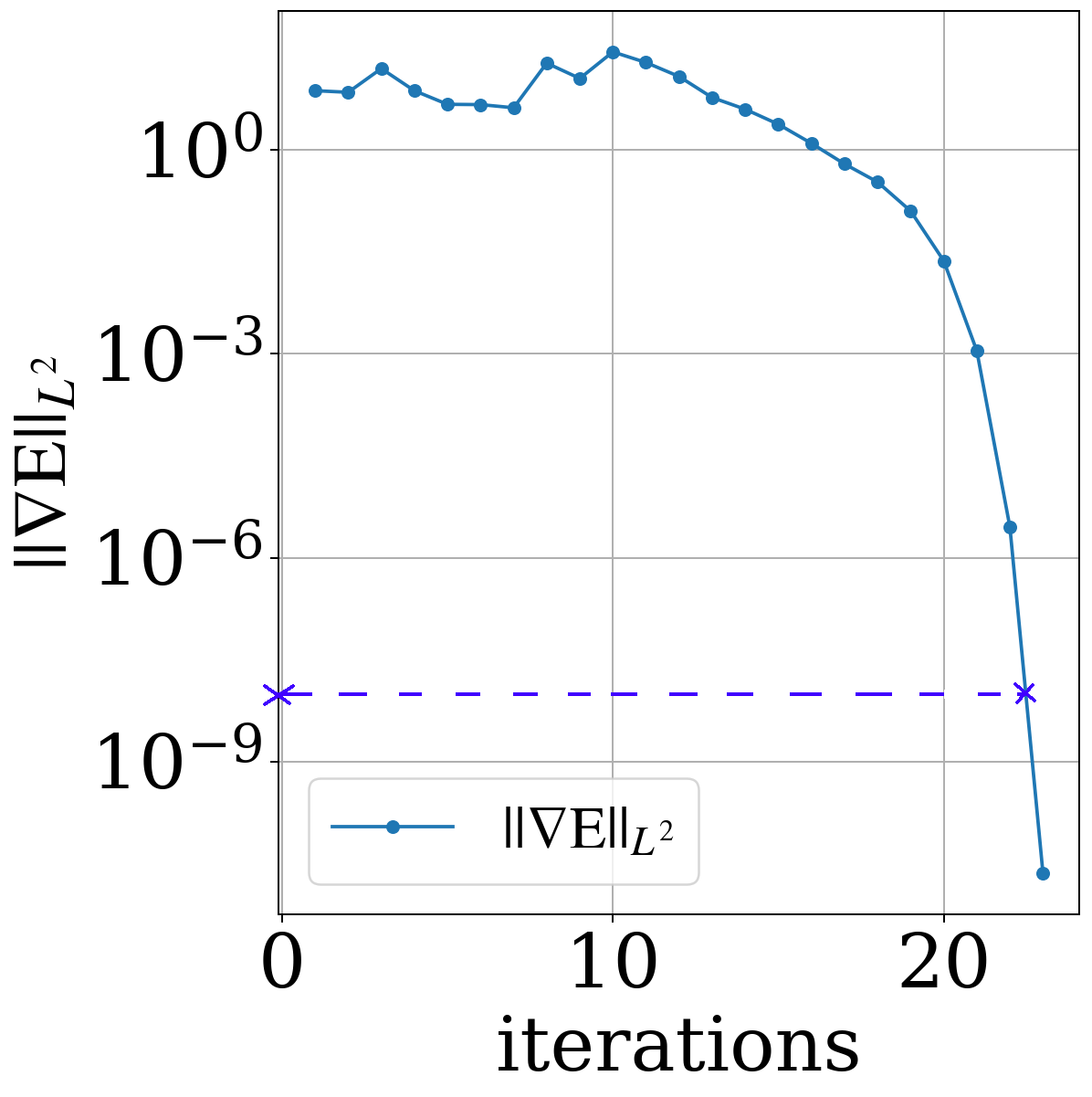}\qquad \includegraphics[height=0.45\linewidth]{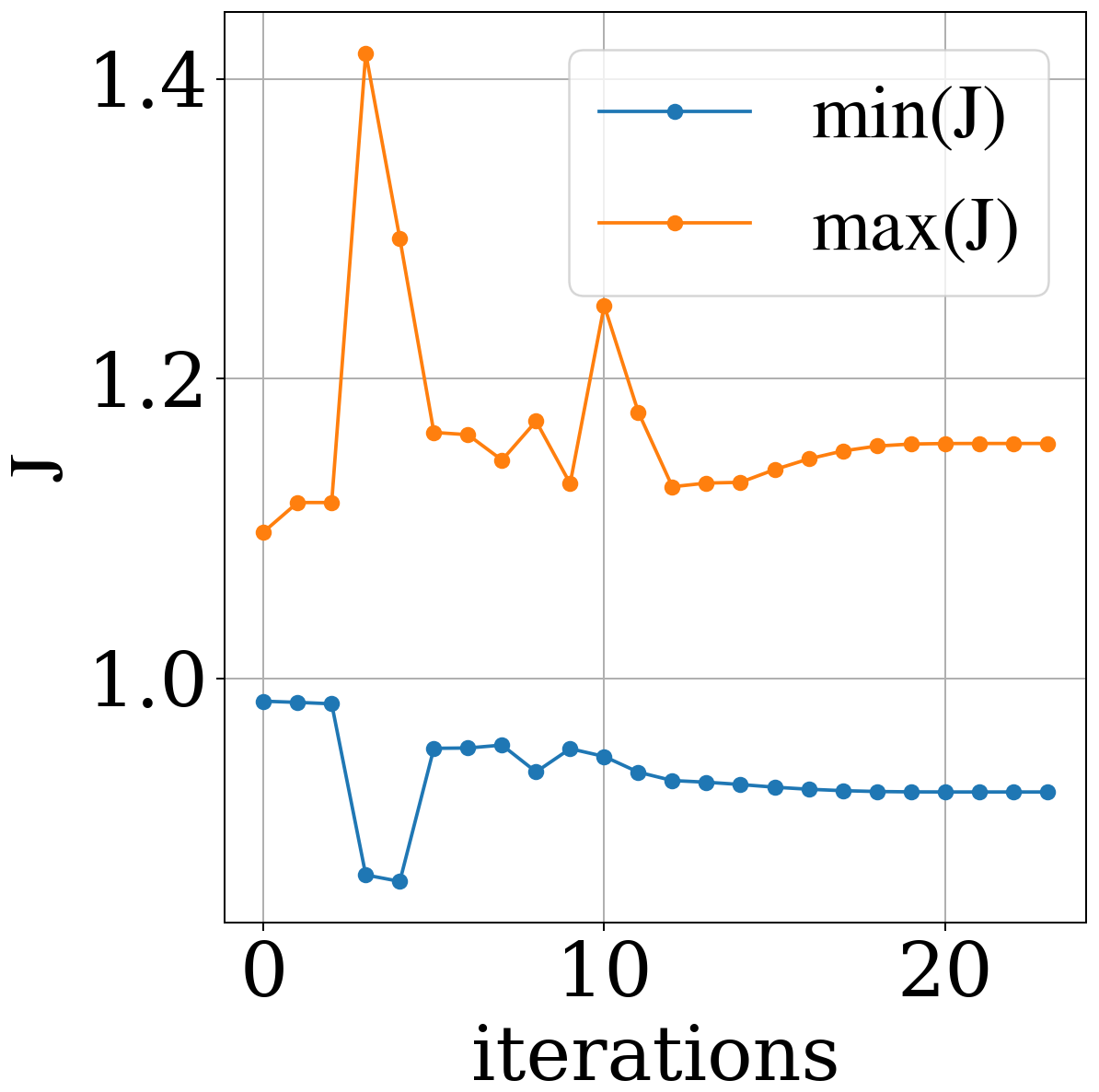}
    \caption{2D forward problem. Left: convergence of the ${\mathrm{L}}^2$ norm of the energy differential. Right: minimum and maximum values of $J$, with respect to Newton iterations.}
    \label{fig:2d-forward-convergence}
\end{figure}

In Figure~\ref{fig:2d-forward-convergence}, the evolution of the ${\mathrm{L}}^2$ norm of the energy differential $\Diff\mathrm{E}(\u;\v)$, together with the minimum and maximum values of the deformation Jacobian, are reported as functions of the iteration number. The results show a progressive decrease in the energy differential norm, reaching $2.27\, 10^{-11}$. Moreover, the deformation Jacobian $J$ remains bounded between $0.85$ and $1.41$, which is consistent with the quasi-incompressibility property.
 
\begin{figure}[htbp!]
    \centering
    \includegraphics[width=\linewidth]{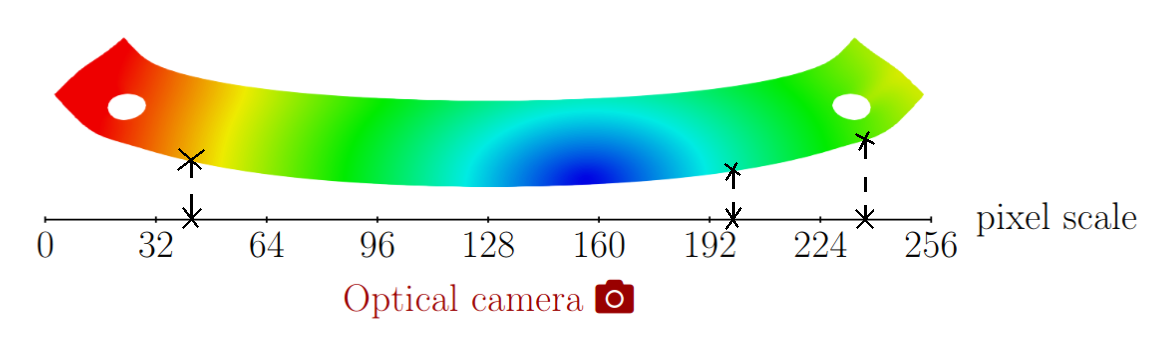}
    \caption{Approximate marker positions chosen on the optical map for the inverse problem with their respective positions on a 1D pixel grid before adding measurement noise. Tissue is colored by the activation time computed using a monodomain model with standard tissue parameters, and fibers from the deformation.}
    \label{fig:2d-inverse-sketch}
\end{figure}

\begin{figure}[htbp!]
    \centering
    \includegraphics[width=0.98\linewidth]{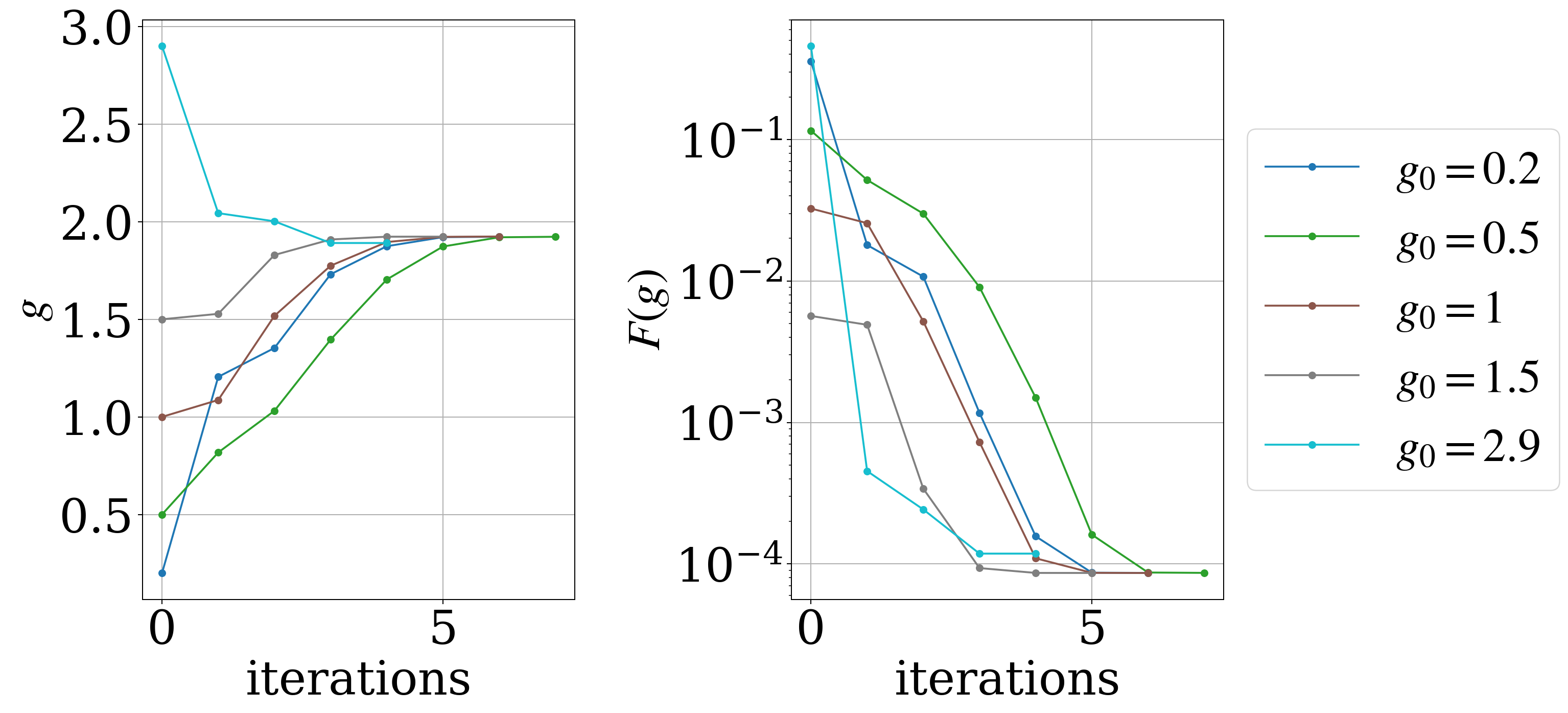}
   \caption{2D inverse problem for several initial guesses of the optimization algorithm. Left: traction amplitude~$g$ with respect to L-BFGS-b iterations. Right: cost function~$F$.}
    \label{fig:2d-inverse-results}
\end{figure}

\textbf{Inverse problem on a 2D domain.} The inverse problem is solved on the same 2D geometry considered in the previous paragraph. As explained in Section~\ref{sec:materials}, the aim is to recover the amplitude of the applied traction from a set of three prescribed marker positions in the optical map ($N_x = 3$, see Figure~\ref{fig:2d-inverse-sketch}). Figure~\ref{fig:2d-inverse-results} reports the evolution of the cost functional given in Equation~\eqref{eq:costFct} and of the  traction amplitude~$g$ during the optimization process. Various initial amplitude guesses are chosen for the L-BFGS-b algorithm. All cases converge toward the same value $g^*=1.92$, close to the reference traction amplitude $g_{\mathrm{ref}}=2$. The difference is explained by the uncertainty $\varepsilon_i$ which was added to the target coordinates.


\begin{figure}[htbp!]
    \centering
    \includegraphics[width=0.67\linewidth]{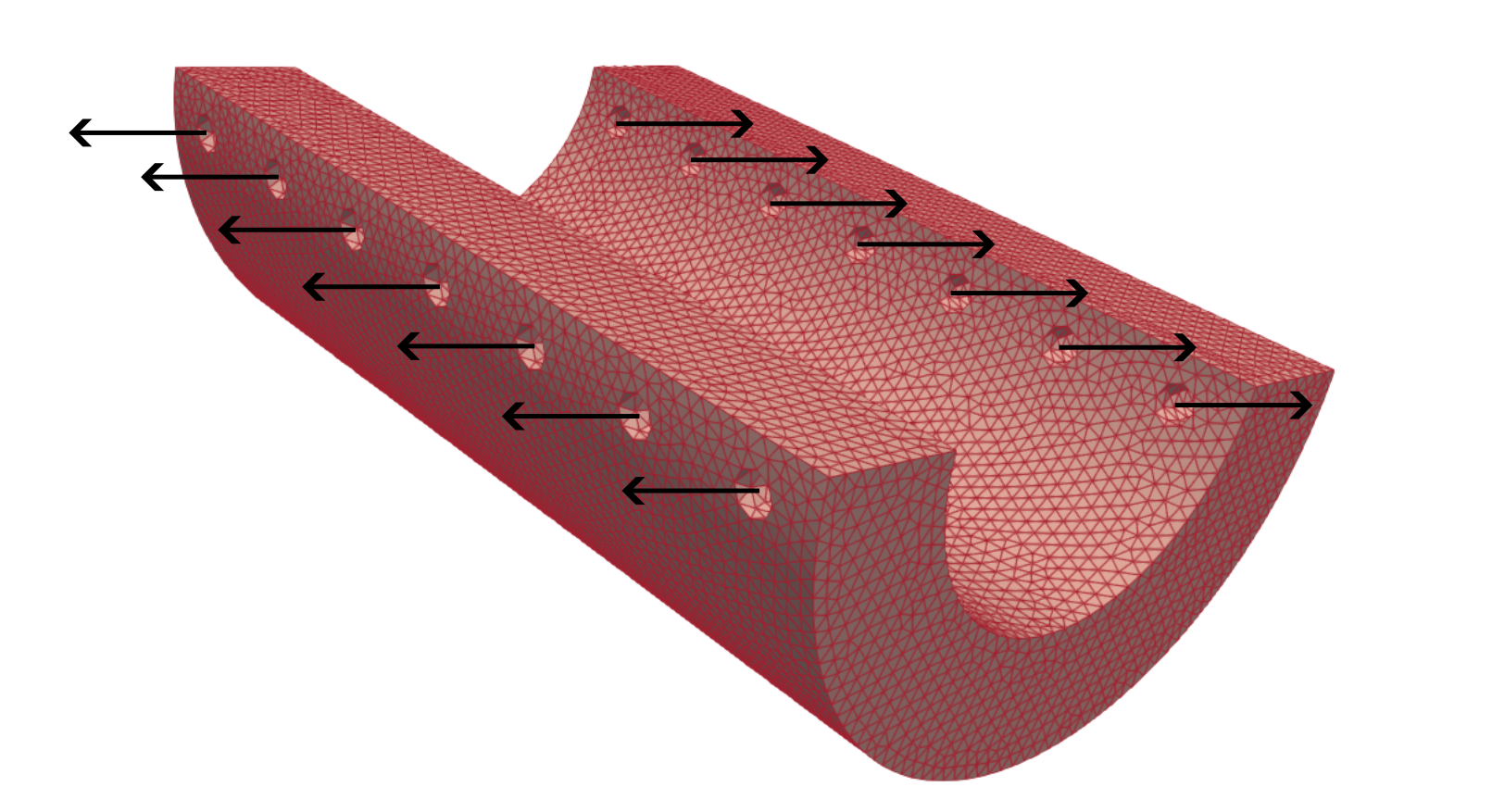} 
    
    \vspace{-0.41cm}
    {\large $\downarrow$}
    \vspace{-0.05cm}
    
    \includegraphics[width=0.64\linewidth]{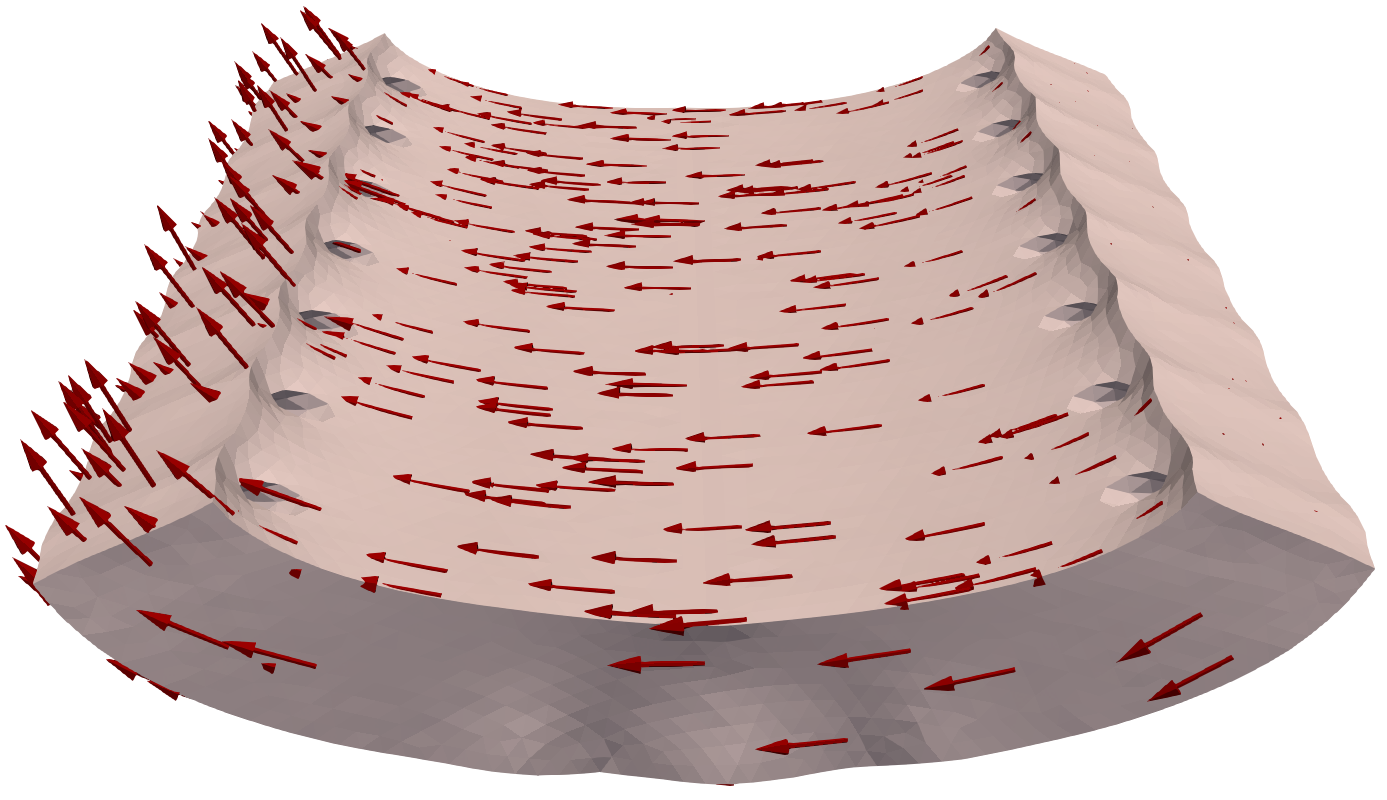} 
    \caption{3D forward problem: reference domain subjected to traction in the $x$-direction and deformed configuration showing a sample of the local fiber directions in red.}
    \label{fig:3d-deformed-fibers}
\end{figure}

\textbf{Extension to a 3D domain.} Here, the domain is a semi-cylindrical shell-like geometry obtained by extruding a lower half-annular along the $z$-direction. Boundary traction forces are applied to the lateral edges of the shell along the $x$-direction, as illustrated in Figure~\ref{fig:3d-deformed-fibers}. The resulting deformed configuration is shown together with a sample of the deformed local fiber directions. The corresponding volume variation is approximately $1.5\%$. The ${\rm L}^2$ norm of the energy differential decreases to $3.46 \, {10^{-8}}$, and the deformation Jacobian remains bounded between $0.87$ and~$1.28$, remaining in the quasi-incompressibility setting.

% \begin{figure}[htbp!]
%     \centering
%     % \includegraphics[width=0.8\linewidth]{Images/3d/energy_Direct_case_10.png}
%     \includegraphics[width=0.8\linewidth]{Images/3d/L2Gradenergy_Direct_case_10.png}
%     \includegraphics[width=0.78\linewidth]{Images/3d/jacobian_Direct_case_10.png}
%     \caption{3D forward problem: evolution, with respect to the iterations, of the norm of the energy differential and of the minimum and maximum values of the Jacobian.}
%     \label{fig:3d-forward-convergence}
% \end{figure}

\section{Discussion and conclusion}

This work presents a biomechanical registration framework to link high-resolution MRI data acquired in a resting configuration to optical-mapping data obtained in a stretched configuration. The proposed approach relies on a large deformation model based on nearly incompressible transversely isotropic hyperelasticity. The 2D forward simulations show that the prescribed boundary traction reproduces a flattened domain with an admissible deformation. 
The inverse problem provides encouraging results, as the optimization converges toward the reference traction amplitude from several initial guesses, using only a limited number of marker positions, consistent with the available data.
The 3D simulation confirms the possibility of extending the method to the right-ventricle geometry in Figure~\ref{fig:3d-rv-mesh}.

\begin{wrapfigure}{r}{0.17\textwidth}
    \centering
    \includegraphics[width=\linewidth]{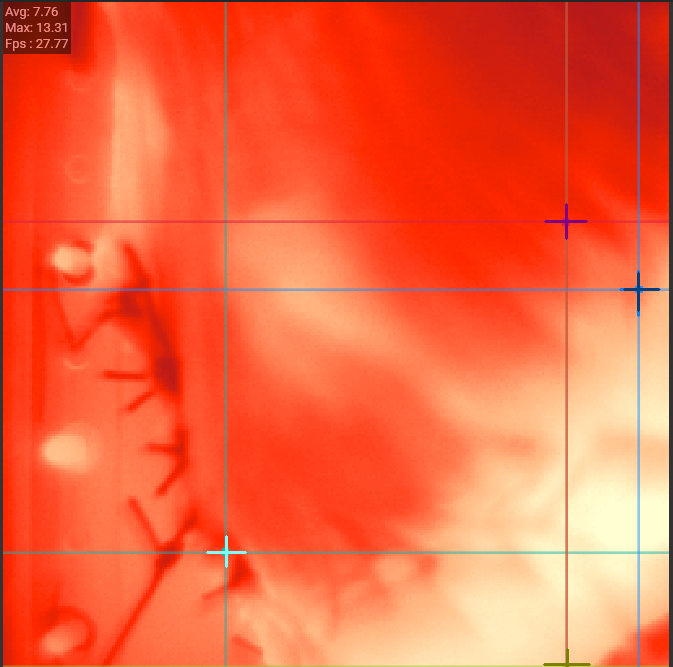}
    \caption{Tissue wedge in the optical mapping setup with position markers}
    \label{fig:optical-map}
\end{wrapfigure}

Several limitations remain. The computational cost is expected to increase substantially in 3D due to both the use of a refined mesh of a realistic geometry (Figure~\ref{fig:3d-rv-mesh}) and the need to solve the forward problem at each iteration of the inverse problem.
Another challenge is the choice of material parameters, adapted from previous studies \cite{Guccione1991,HolzapfelOgden2009} and adjusted to the present framework. Their calibration remains delicate and requires further validation.

\textcolor{black}{This paper addressed a simple 2D inverse problem involving a continuous functional that is likely convex. Future work will focus on extending the method to the more complex right-ventricle geometry in Figure~\ref{fig:3d-rv-mesh} and on using the limited number of experimental marker positions extracted from optical-mapping images (Figure~\ref{fig:optical-map}). Moreover, the inverse problem will no longer be limited to estimating the force amplitude, but will also require identifying the positions and directions of the forces.}

\section*{Acknowledgments}
% This work was supported by the H2020 EU SimCar-dioTest project (Digital transformation in Health and Care SC1-DTH-06-2020; grant agreement number 101016496). 
This study was supported by a research grant from the French Government, managed
by the National Funding Agency for Research (ANR), 10-IAHU-04. The computations are provided by the computing facilities of the MCIA (Mésocentre de Calcul Intensif Aquitaine). We would also like to thank Richard Walton for insightful and constructive discussions. 

% LateX generates the ``References'' heading automatically and switches
% to 9 point type for the bibliography.  Please  use BibTeX and
% follow the examples in the sample 'refs.bib' file to enter your references.
\bibliography{refs_simple}

\begin{thebibliography}{1}
\expandafter\ifx\csname url\endcsname\relax
  \def\url#1{\texttt{#1}}\fi
\expandafter\ifx\csname urlprefix\endcsname\relax\def\urlprefix{URL }\fi

\bibitem{Bone2018Deformetrica}
B{\^o}ne A, {\it et. al. }.
\newblock Deformetrica 4: An open-source software for statistical shape
  analysis.
\newblock In Shape Med Imaging, volume 11167 of Lect. Notes Comput. Sci. Cham:
  Springer, 2018;\hspace{0pt} 3--13.

\bibitem{Genet2018EquilibratedWarping}
Genet M, {\it et. al. }.
\newblock Equilibrated warping: Finite element image registration with finite
  strain equilibrium gap regularization.
\newblock Med Image Anal 2018;\hspace{0pt}50:1--22.

\bibitem{walton2018compartmentalized}
Walton RD, {\it et. al. }.
\newblock Compartmentalized structure of the moderator band provides a unique
  substrate for macroreentrant ventricular tachycardia.
\newblock Circ Arrhythm Electrophysiol 2018;\hspace{0pt}11(8):e005913.

\bibitem{pannetier2023modeling}
Pannetier V, {\it et. al. }.
\newblock Modeling cardiac stimulation by a pacemaker, with accurate
  tissue-electrode interface.
\newblock In FIMH, volume 13958 of Lect. Notes Comput. Sci. Springer Nature
  Switzerland, 2023;\hspace{0pt} 194--203.

\bibitem{HolzapfelOgden2009}
Holzapfel GA, Ogden RW.
\newblock Constitutive modelling of passive myocardium: A structurally based
  framework for material characterization.
\newblock Philos Trans R Soc A 2009;\hspace{0pt}367(1902):3445--3475.

\bibitem{gasser2006hyperelastic}
Gasser TC, {\it et. al. }.
\newblock Hyperelastic modelling of arterial layers with distributed collagen
  fibre orientations.
\newblock J R Soc Interface 2006;\hspace{0pt}3(6):15--35.

\bibitem{Guccione1991}
Guccione JM, {\it et. al. }.
\newblock Passive material properties of intact ventricular myocardium
  determined from a cylindrical model.
\newblock J Biomech Eng 1991;\hspace{0pt}113(1):42--55.

\bibitem{ogden1978nearly}
Ogden RW.
\newblock Nearly isochoric elastic deformations: Application to rubberlike
  solids.
\newblock J Mech Phys Solids 1978;\hspace{0pt}26(1):37--57.

\bibitem{Fung1993Biomechanics}
Fung YC.
\newblock Biomechanics: Mechanical Properties of Living Tissues.
\newblock 2 edition. New York: Springer-Verlag, 1993.

\end{thebibliography}

\begin{correspondence}
Joyce Ghantous\\
Univ. Bordeaux, CNRS, Inria, Bordeaux INP, IMB, UMR 5251, IHU Liryc, F-33400 Talence, France\\ % full postal adresse
joyce.ghantous@inria.fr
\end{correspondence}

\end{document}